\documentclass[12pt]{amsart}
\usepackage{graphicx} % Required for inserting images
\usepackage{booktabs}
\usepackage{float}
\usepackage{soul}
\usepackage[citecolor=Navy]{hyperref}
\usepackage[ruled,vlined]{algorithm2e}
\usepackage{amsfonts,amssymb,amsmath,amsthm}
\usepackage{algorithmic}
\usepackage{amssymb,amscd,amsxtra,calc}
\usepackage{cmmib57}
\usepackage{graphicx}
\usepackage{mathrsfs,enumerate,xcolor}

\usepackage[all]{xy}
\usepackage{multirow} % Required for multirows
\usepackage{psfrag,xmpmulti,amscd,color,pstricks, import}
\usepackage{tikz-cd} 

\newtheorem{thm}{Theorem}[section]

\theoremstyle{definition}

\newtheorem{example}[thm]{Example}
\newtheorem*{rem*}{Remark}

\title[BNQN, Schr\"oder's theorem,  and Linear Conjugacy]{Backtracking New Q-Newton's method, Schr\"oder's theorem,  and Linear Conjugacy}
\date{\today}

\author[Forn\ae ss]{John Erik Forn\ae ss}
\address{Department of Mathematics, NTNU, Norway }
\email{fornaess@gmail.com}

\author[Hu]{Mi Hu}
\address{Department of Mathematics, University of Oslo, Norway}
\email{humihqu@gmail.com}

\author[Truong]{Tuyen Trung Truong}
\address{Department of Mathematics, University of Oslo, Norway}
\email{tuyentt@math.uio.no}

\author[Watanabe]{Takayuki Watanabe}
\address{Chubu University Academy of Emerging Sciences/Center for Mathematical Science and Artificial Intelligence, Chubu University, Japan}
\email{takawatanabe@isc.chubu.ac.jp}

\begin{document}

\maketitle

\begin{abstract}
A new variant of Newton's method - named Backtracking New Q-Newton's method (BNQN) - which has strong theoretical guarantee, is easy to implement, and has good experimental performance, was recently introduced by the third author. 

Experiments performed previously showed some remarkable properties of the basins of attractions for finding roots of polynomials and meromorphic functions using BNQN. In particular, it seems that for finding roots of polynomials of degree 2, the basins of attraction of the dynamics for BNQN are the same as that for Newton's method (the latter is the classical Schr\"oder's result in Complex Dynamics). 

In this paper, we show that indeed the picture we obtain when finding roots of polynomials of degree 2 is the same as that in Sch\"oder's result, with a remarkable difference: on the boundary line of the basins, the dynamics of Newton's method is chaotic, while the dynamics of BNQN is more smooth. On the way to proving the result, we show that BNQN (in any dimension) is invariant under conjugation by linear operators of the form $A=cR$, where $R$ is unitary and $c>0$ a constant. This again illustrates the similarity-difference relation between BNQN and Newton's method.  

\end{abstract}

\section{Introduction}

Effectively finding roots of polynomials is an extensive topic within many fields, e.g. Mathematics,  Optimization and Computer Science. In particular, it was one origin of Complex Dynamics and is still extensively researched nowadays.   

It is known since Abel and Ruffini's time that for polynomials of degree at least 5 in one complex variable, there is no closed form formulars in terms of radicals of the coefficients to find their roots. There are closed form formulas for the roots, but they are in terms of transcendental functions \cite{RefU}, and hence can be difficult to use in practical situations.  This means that the most one can do is to find approximate solutions, and must rely on iterative methods. 

A well known iterative method for finding roots of polynomial equations $P(z)=0$ is Newton's method, which works as follows. One chooses a random initial point $z_0\in \mathbf{C}$, and then defines a sequence 
\begin{eqnarray*}
z_{n+1}=z_n-\frac{P(z_n)}{P'(z_n)}. 
\end{eqnarray*}
One attractive property of Newton's method is that it has a fast local rate of convergence: if $z_0$ is  close enough to a simple root $z^*$ of $P(z)$, then the constructed sequence $\{z_n\}$ will converge to $z^*$, and the rate of convergence is quadratic. 

However, Newton's method has no guarantee for global convergence: if the initial point $z_0$ is not close to a root of $P(z)$, then $\{z_n\}$ may not converge to any root of $P(z)$. Indeed, Newton's method applied to a polynomial in one complex variable gives rise to an algebraic dynamical system on the Riemann sphere, and McMullen \cite{RefMc} showed that no algebraic dynamical system can find roots of a general polynomial of degree at least 4. A consequence of this result is that even to find roots of only polynomials, one needs to rely on non-algebraic methods. 

Another issue of Newton's method is its sensitivity to the initial point. This is evident by the fact that the basins of attraction for Newton's method are usually fractal \cite{RefN}. This fact has both good and bad consequences. A good consequence is that it produces beautiful pictures. A bad consequence is that it is then difficult to make an intelligent prediction on where the constructed sequence will end up, given an approximate information of the initial point. 

This paper concerns the dynamics of iterative methods which are close variants of Newton's method and which have both strong convergence guarantees for finding roots polynomial or meromorphic functions in 1 complex variable, and easy to implement on computers. As such, and to keep the paper concise, we apologize for not being able to discuss/review the many thousands variations of Newton's method in the current literature. However, readers who are interested in such reviews and histories, can easily find plenty references (e.g. in \cite{{RefAlex}}, \cite{{RefTT}} and \cite{RefT}).    

A direct generalization of Newton's method is Relaxed Newton's method. This method works as follows. One chooses a non-zero complex number $\alpha $ and construct a sequence 
\begin{eqnarray*}
z_{n+1}=z_n-\alpha \frac{P(z_n)}{P'(z_n)}. 
\end{eqnarray*}
Relaxed Newton's method, which is still an algebraic dynamical system when applied to a polynomial, again cannot have global convergence guarantee, according to the cited result \cite{RefMc}. A surprising result, by Sumi \cite{RefS}, is that if one adds {\bf randomness} into the design of Relaxed Newton's method (i.e. if one chooses $\alpha$ randomly in a precise manner) then one has global convergence when applied to polynomials. The precise statement of the result is recalled in Section 2.    
It is unknown if global convergence of Random Relaxed Newton's method also holds when used for transcendental holomorphic functions or more generally meromorphic functions in 1 complex variable. Another fact is that the rate of convergence may not be quadratic anymore, given that the $\alpha _n$ must be uniformly randomly chosen from a set far from $1$. Moreover, very little is known about its guarantees in higher dimensions or for finding roots in real numbers.  

\begin{figure}
\centering
%\includegraphics[scale=.5]{BasinOfAttractionDegree2BacktrackingGD.png}\hfill
%\caption{For Newton's method}\hfill
 \includegraphics[scale=.3]{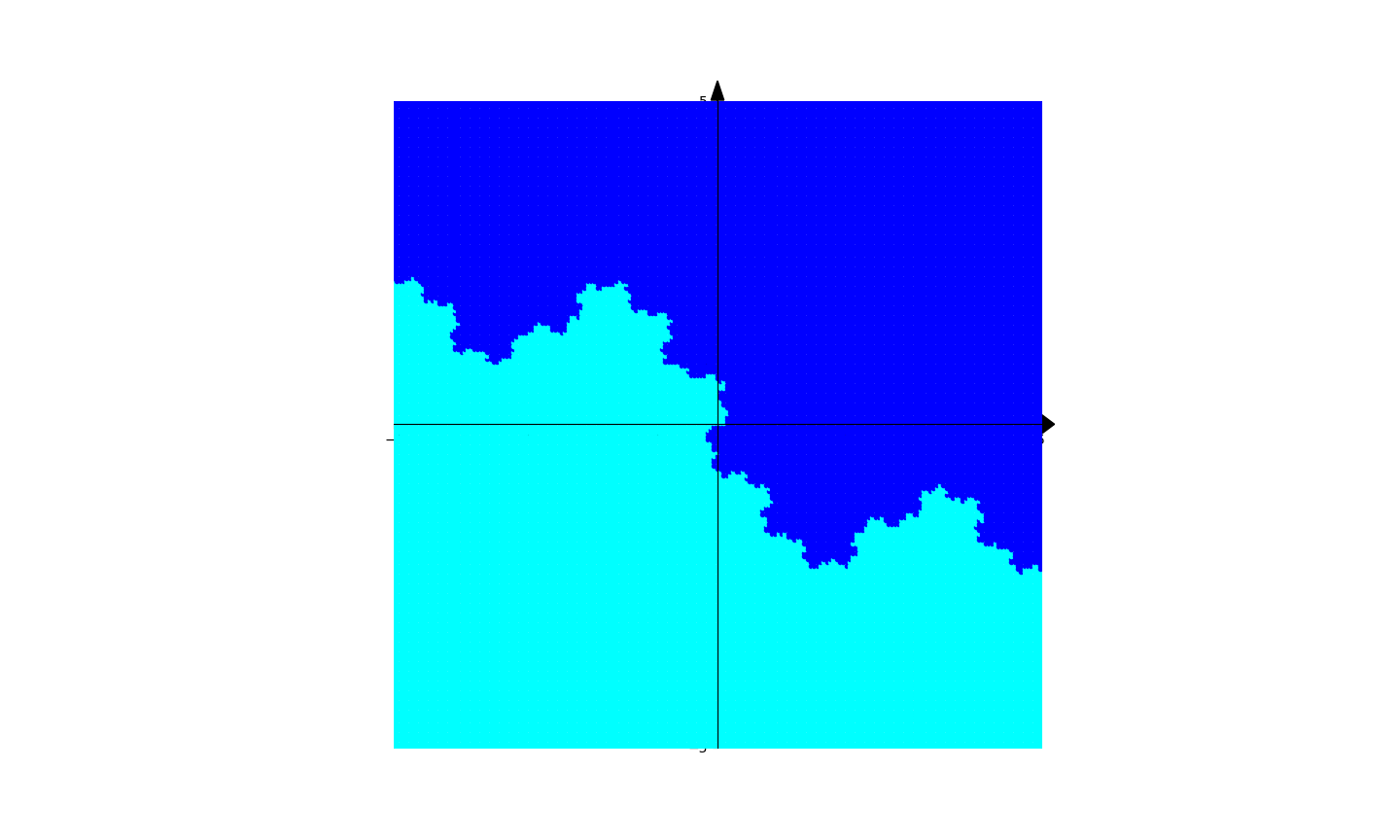}%\hfill
\caption{Basins of attraction for finding roots of  a polynomial  of degree 2, for one run of  Random Relaxed Newton's method.  Initial points of the same color produces sequences converging to the same root of the polynomial. Different runs of Random Relaxed Newton's method produce different pictures.}
 \label{fig:Degree2RandomRelaxedNewtonMethod} 
\end{figure}

\begin{figure}
\centering
%\includegraphics[scale=.5]{BasinOfAttractionDegree2BacktrackingGD.png}\hfill
%\caption{For Newton's method}\hfill
% \includegraphics[scale=.3]{BasinOfAttractionDegree2RelaxingNewtonMethodFixedFactor25August2023.png}\hfill
 \includegraphics[scale=.5]{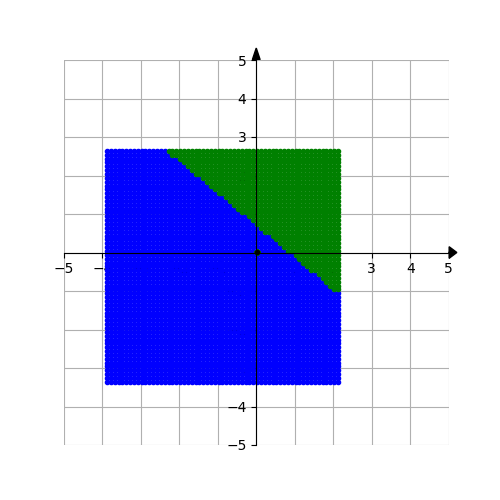}
\caption{Basins of attraction for finding roots of a polynomial  of degree 2. BNQN produces the same picture as Newton's method. Initial points of the same color produces sequences converging to the same root of the polynomial. The behaviours of the two methods on the boundary line are on the other hand quite different, as discussed in detail in this paper's text.}
 \label{fig:Degree2} 
\end{figure}

Similar to Newton's method, the basins of attraction for Random Relaxed Newton's method are usually also very fractal. Moreover, the basins of attraction depend on the choice of the sequence $\{\alpha _n\}$, which makes it difficult to talk about a unique basin of attraction for this method. Basins of attraction for this method can be very chaotic, and it seems that they are distortions of basins of attraction of the vannila Newton's method. Even for a polynomial of degree 2, the basin of attraction can be drastically different from the familiar picture in Schr\"oder's theorem \cite{RefSE2} (and see Section 2 for a precise statement of this theorem), see Figures \ref{fig:Degree2RandomRelaxedNewtonMethod} and \ref{fig:Degree2}. The fourth author (in ongoing work) has been carrying many experiments concerning basins of attraction for Random Relaxed Newton's method.  

With the above brief review on Newton's method, we are now ready to discuss the main topic of this paper: Backtracking New Q-Newton's method (BNQN) \cite{RefT}. The method has the best theoretical guarantee so far, among variants of Newton's method, concerning finding roots of meromorphic functions in one complex variable.  This method can be applied to the more general setting of solving optimization problems in any dimension. The precise definition of the algorithm and some of its relevant main properties are given in Section 2. In \cite{RefT}, many experimental results are presented, which show a remarkable phenomenon: it seems that in general the basins of attraction for BNQN are more smooth than that of Newton's method, and seem to reflect the geometric configuration of the roots of the concerned function. This leads naturally to the question of whether what observed in those experiments is only a lucky coincidence of numerical effects, or indeed reflects some very deep reasons. 

This paper is the first step to providing rigorously theoretical justifications to these experimental results. We treat the first non-trivial case, that of finding roots of a polynomial of degree 2. In \cite{RefT}, there are two distinct versions of BNQN depending on whether the objective function has compact sublevels or not (see Section 2). In this paper, we introduce a new variant, named BNQN New Variant, which includes both these versions as special case, see Section 2 for detail. Our first main result is the following. 

\begin{thm} Let $g(z)$ be a polynomial of degree $2$. Let $F:\mathbf{R}^2\rightarrow \mathbf{R}$ be defined as $F(x,y)=|g(x+iy)|^2/2$. 

For an initial point $z_0\in \mathbf{C}$, let $\{z_n\}$ be the sequence constructed by BNQN New Variant for the cost function $F$. 

1) Assume that $g(z)=c(z-z^*)^2$, where $c,z^*\in \mathbf{C}$. Then, for every initial point $z_0\in \mathbf{C}$, the constructed sequence $\{z_n\}$ converges to the (multiplicity 2) root $z^*$ of $g$. 

2) Assume that $g(z)=c(z-z_1^*)(z-z_2^*)$, where $c,z_1^*,z_2^*\in \mathbf{C}$ and $z_1^*\not= z_2^*$. Let $L$ be the perpendicular bisector of the interval connecting $z_1^*$ and $z_2^*$, and $H_1,H_2\subset \mathbf{C}$ be two half-planes divided by $L$, where $z_1^*\in H_1$ and $z_2^*\in H_2$.  

a) If the initial point $z_0$ is in $H_1$, then the constructed sequence $\{z_n\}$ converges to $z_1^*$. Similarly, if the initial point $z_0$ is in $H_2$, then the constructed sequence $\{z_n\}$ converges to $z_2^*$.

b) On the other hand, if the initial point $z_0$ is in $L$, then the constructed sequence $\{z_n\}$ converges to $z^*=(z_1^*+z_2^*)/2$, where $z^*$ is the critical point of $g$, that is $g'(z^*)=0$.  
\label{TheoremMain1}\end{thm}

The conclusion of part 2a of Theorem \ref{TheoremMain1} is the same as that of Schr\"oder's theorem for the dynamics of Newton's method for a polynomial of degree 2. On the other hand, the conclusion of part 2b of Theorem \ref{TheoremMain1} is different from that of Newton's method, which is chaotic on the line $L$, see Section 2 for detail.   

The proof of Theorem \ref{TheoremMain1} uses the fact that BNQN (New Variant) is invariant under linear conjugations by rotations. This allows us to reduce to considering only two polynomials: $f(z)=z^2$ and $f(z)=z^2-1$, much like the strategy in the proof of Schr\"oder's theorem. This invariant property can be extended to a more general form in higher dimension and is itself of independent interest. 

\begin{thm} The dynamics of BNQN New Variant is invariant under conjugation by linear operators of a certain form. The precise statement is as follows. 

Let $F:\mathbf{R}^m\rightarrow \mathbf{R}$ be an objective function. Let $A:\mathbf{R}^m\rightarrow \mathbf{R}^m$ be an invertible linear map, of the form $A=cR$, where $R$ is  {\bf unitary} (i.e. $RR^T=Id$, where $R^T$ is the transpose of $R$) and $c>0$. Define $G(z)=F(Az)$. 

Let $z_0\in \mathbf{R}^m$ be an initial point, and let $\{z_n\}$ be the sequence constructed by BNQN (with parameters $\delta _0,\delta _1,\ldots ,\delta _m, \theta, \tau$) applied to the function $F(z)$. 

Let $z_0'=A^{-1}z_0$, and let $\{z_n'\}$ be the sequence constructed by BNQN (with parameters $\delta _0',\delta _1',\ldots ,\delta _m', \theta ', \tau $) applied to the function $G(z)$. 

Assume that $\delta_i'=\delta_i c^{2-\tau}$ for all $0\leq i\leq m$, and $\theta '=c\theta $.

Then for all $n$, we have $z_n'=A^{-1}z_n.$

\label{TheoremMain2}
\end{thm}

We recall that Newton's method is invariant under every invertible linear change of coordinates, see Section 2 for detail.  After the proof of Theorem \ref{TheoremMain2}, we will present a simple example showing that without the condition on the linear operator $A$, then the conclusion of the theorem may not hold. Therefore, here we see again that there is a similarity-difference relation between the dynamics of BNQN and that of Newton's method.

{\bf Remark:} 

While in general BNQN (New Variant) needs that the parameters $\delta _0,\ldots ,\delta _m$ in its definition are random to ensure avoidance of saddle points, the proofs of Theorems \ref{TheoremMain1} and \ref{TheoremMain2} do not need that. It is also interesting to see the appearance of negative inverse of the golden ratio $(\sqrt{5}+1)/2$, i.e. the number $(1-\sqrt{5})/2$, which is crucially used at some places in the proof of Theorem \ref{TheoremMain1}.

\begin{figure}
%\centering
\includegraphics[scale=.3]{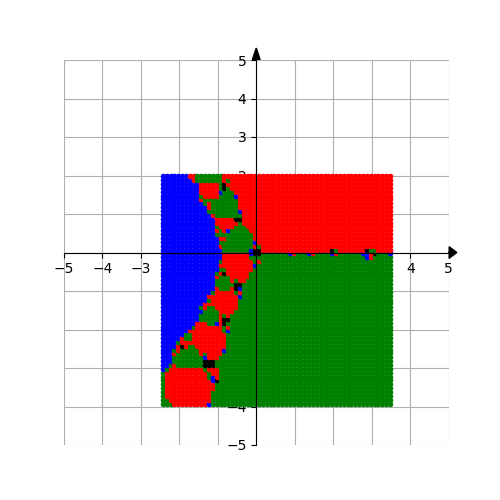}\hfill
%\caption{For Newton's method}\hfill
 \includegraphics[scale=.3]{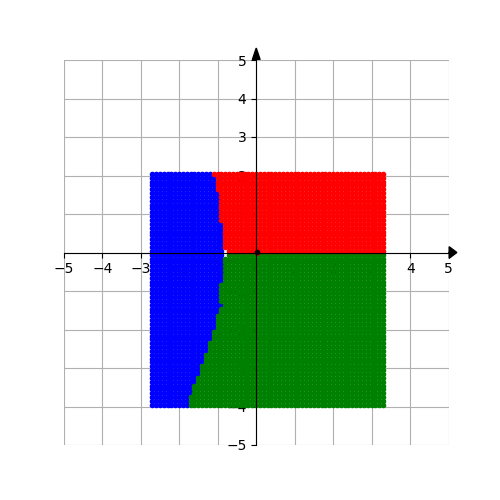}\hfill
 \includegraphics[scale=.3]{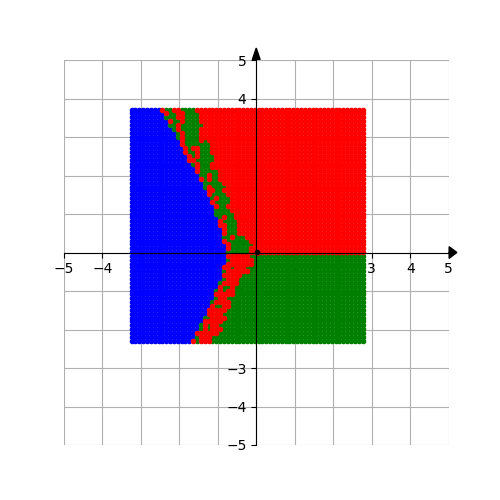}
\caption{Basins of attraction for finding roots of a polynomial  of degree 3. The left image is for Newton's method, the central image is for BNQN, the right image is for Backtracking line search for Gradient Descent.  Initial points of the same color produces sequences converging to the same root of the polynomial.}
 \label{fig:Degree3} 
\end{figure}

The proof of Theorem \ref{TheoremMain1} can be modified to show that when applied to finding roots of a polynomial of degree 2, Backtracking line search for Gradient Descent (see \cite{RefAr}) also has the same basins of attraction as in BNQN. However, things are different for polynomials of degree 3 and higher, see Figure \ref{fig:Degree3}. It seems to us that one reason for the difference lies in the fact that in BNQN one uses the Hessian matrix, which usually carry additional information about the curvature of the underlying function landscape. 

This is the first part, concentrating more on the theoretical part, of BNQN. A following up paper, concentrating more on the experimental part, will explore in depth the relations between BNQN and some other more established topics in ODE and geometry. 

{\bf Plan of the paper:} In Section 2 we recall some relevant facts on Newton's method, Random Relaxed  Newton's method and BNQN. In Section 3 we provide the proofs of Theorems \ref{TheoremMain1} and \ref{TheoremMain2}, together with some examples and further comments. 

{\bf Acknowledgments:} Mi Hu and Tuyen Trung Truong are partially supported by Young Research Talents grant 300814 from Research Council of Norway. Takayuki Watanabe is partially supported by JSPS Grant-in-Aid for Early-Career Scientists Grant Number JP 23K13000. Some ideas of the work were discussed/carried out in the inspiring environments of the conferences ICIAM 2023 (Waseda University) and ``Birational geometry and algebraic dynamics" 2023 (University of Tokyo), for which we would like to thank the organisers and participants, in particularly Hiroki Sumi, Mark Comerford and Fabrizio Catanese for helpful comments. 

\section{Preliminaries} In this section, we briefly review properties of Newton's method, Random Relaxed Newton's method, and BNQN.  In particular,  properties of BNQN which will be used later will be presented in enough detail. 

\subsection{Newton's method, Linear conjugacy invariance and Schr\"oder theorem}

We recall here two relevant properties of Newton's method: the first is its invariance under linear change of coordinates, and the second is Schr\"oder's theorem when using it to find roots of polynomials of degree 2. 

We will present the linear transformation invariance of Newton's method in the optimization setting also, to easily compare later to BNQN. We first recall this algorithm. 

Assume that one wants to find (local) minima of an objective function $F:\mathbf{R}^m\rightarrow \mathbf{R}$. Let $\nabla F$ be the gradient of $F$, and $\nabla ^2F$ be the Hessian of $F$. Then Newton's method for the function $F$ is the following iterative algorithm: Choose $z_0\in \mathbf{R}^m$ an initial point, and define: 
$$z_{n+1}=z_n-(\nabla ^2F(z_n))^{-1}.\nabla F(z_n).$$

\begin{thm} Newton's method (both the version for finding root of a complex function $f$ in 1 complex variable and the version for optimizing an objective real function $F$) is invariant under a linear change of coordinates. More precisely, 

1) Let $f(z):\mathbf{C}\rightarrow \mathbf{C}\cup \{\infty\}$ be a meromorphic function in 1 complex variable. Let $a\in \mathbf{C}^*$ be a constant and consider $g(z)=f(az)$. 

Let $z_0\in \mathbf{C}$ be an initial point, and let $\{z_n\}$ be the sequence constructed by Newton's method applied to the function $f(z)$. 

Let $z_0'=z_0/a$, and let $\{z_n'\}$ be the sequence constructed by Newton's method applied to the function $g(z)$. 

Then for all $n$, we have $z_n'=z_n/a$. 

2) Let $F:\mathbf{R}^m\rightarrow \mathbf{R}$ be an objective function. Let $A:\mathbf{R}^m\rightarrow \mathbf{R}^m$ be an invertible linear map. Define $G(z)=F(Az)$. 

Let $z_0\in \mathbf{R}^m$ be an initial point, and let $\{z_n\}$ be the sequence constructed by Newton's method applied to the function $F(z)$. 

Let $z_0'=A^{-1}.z_0$, and let $\{z_n'\}$ be the sequence constructed by Newton's method applied to the function $G(z)$. 

Then for all $n$, we have $z_n'=A^{-1}.z_n$. 
\label{TheoremNewtonLinearConjugacy}\end{thm}
\begin{proof} The proof is standard, however we present here for the readers' convenience, and also for to compare with the proof used in the case of BNQN. See \cite{RefM} for more detail.  

1) See e.g. \cite{RefM}. By definition of Newton's method, we have
\begin{eqnarray*}
 z_{n+1}&=&z_n-\frac{f(z_n)}{f'(z_n)},\\
 z_{n+1}'&=&z_n'-\frac{g(z_n')}{g'(z_n')}. 
\end{eqnarray*}

By definition of $g$, we have $g'(z)=a.f'(az)$. From this, it is easy to verify that $z_n'=z_n/a$ for all $n$.

2) Let $A^T$ be the transpose of $A$. Then we can check that 
\begin{eqnarray*}
 \nabla G(z)&=&A^T.\nabla F(Az),\\
 \nabla ^2G(z)&=&A^T.\nabla ^2F(Az).A.
\end{eqnarray*}
By definition of Newton's method, we have
\begin{eqnarray*}
 z_{n+1}&=&z_n-(\nabla ^2F(z_n))^{-1}.\nabla F(z_n),\\ 
 z_{n+1}'&=&z_n'-(\nabla ^2G(z_n'))^{-1}.\nabla G(z_n')\\
 &=&z_n'-(A^{-1}.(\nabla ^2 F(Az_n'))^{-1}.(A^T)^{-1}).(A^T.\nabla F(Az_n'))\\
 &=&z_n'-A^{-1}.(\nabla ^2 F(Az_n'))^{-1}.\nabla F(Az_n')\\
 &=&A^{-1}.(Az_n'-(\nabla ^2 F(Az_n'))^{-1}.\nabla F(Az_n')). 
\end{eqnarray*}
Hence, if $z_0'=A^{-1}z_0$, then we can show by induction that $z_{n}'=A^{-1}z_n$ for all $n$. 
    
\end{proof}

We next explore the behaviour of Newton's method when applied to finding roots of a polynomial of degree 2. 

\begin{thm} Let $f(z)$ be  a polynomial of degree $2$. Let $z_0\in \mathbf{C}$ be an initial point, and let $\{z_n\}$ be the sequence constructed by Newton's method applied to finding roots of $f(z)$. 

1) If $f(z)=c(z-z^*)^2$, where $c\in \mathbf{C}$ is a constant, then for every $z_0$, the sequence $\{z_n\}$ converges to the root $z^*$ of $f(z)$.

2) {\bf Schr\"{o}der's theorem}\cite{RefSE2}  Let $f(z)=c(z-z_1^*)(z-z_2^*)$, where $c,z_1^*,z_2^*\in \mathbf{C}$ and $z_1^*\not= z_2^*$. L $L$ be the 
 perpendicular bisector of the line segment connecting $z_1^*$ and $z_2^*$. Let $H_1,H_2\subset\mathbf{C}$ be the open half-planes separated by $L$, where $z_1^*\in H_1$ and $z_2^*\in H_2$. 

a) If $z_0\in H_1$, then the constructed sequence $\{z_n\}$ converges to $z_1^*$.  Similarly, if $z_0\in H_2$, then the constructed sequence $\{z_n\}$ converges to $z_2^*$. 

b) On the other hand, the behaviour of Newton's method on the line $L$ is very chaotic. In particular: 

- There are finite orbits of any length, and the set of periodic points is dense.  

- There are points whose orbit is dense in $L$.

\label{TheoremSchroderNewton}\end{thm}

\begin{proof} We note that Newton's method is invariant also under translations $z\mapsto z+b$. Also, Newton's method applied to a map $cf(z)$ is the same as Newton's method applied to $f(z)$, for all $c\in \mathbf{C}\backslash \{0\}$. This, and the invariance of Newton's method under linear change of coordinates (part a of Theorem \ref{TheoremNewtonLinearConjugacy}) allows to reduce to considering two special functions $f(z)=z^2$ and $f(z)=z^2-1$. 

1) In this case, as discussed above, we can assume that $f(z)=z^2$. The update rule for Newton's method is then: $z_{n+1}=z_n/2$. Thus, for every initial point $z_0$, the constructed sequence $\{z_n\}$ converges to $0$, which is the unique root of $f(z)$. 

2)

a) Consider $
w=\phi(z)=\frac{z-1}{z+1}
$, which is an automorphism of $\mathbf{P}^1=\mathbf{C}\cup \{\infty\}$, with inverse  $z=\psi(w)=\frac{1+w}{1-w}$. This conjugates the function $f(z)=z^2-1$ to the function 
$$
h(w):=\phi\circ N\circ \psi(w)
     =\phi \circ N\circ \frac{1+w}{1-w}
     =\phi \bigg(\frac{1+w}{2(1-w)}+\frac{1-w}{2(1+w)}\bigg)=w^2.
$$
Hence, the dynamics of $f(z)$ is isomorphic to the dynamics of $h(w)$. Note that $\mathbf{P}^1$ is divided into two connected open sets by the unit circle $\{w:~|w|=1\}$. The set $\{w:~|w|<1\}$ is the basin of attraction for the root $w=0$ of $h(w)$. The set $\{w:~|w|>1\}$ is the basin of attraction for the root $w=\infty$ of $h(w)$. Note that the unit circle $\{w:~|w|=1\}$ is the image of the line $\{z:~Re(z)=0\}$. Therefore, translating back to the z coordinates, we obatin the claims of part a) of the theorem.

b) This is classically known. For details the readers can see for example Example 8.6 in the book \cite{RefD}.  

\end{proof}

\subsection{Random relaxed Newton's method} We recall here the relevant result in \cite{RefS} on convergence of Random Relaxed Newton's method when applied to finding roots of polynomials. 

\begin{thm} Let $0.5<\rho <1$ be a constant. Let $P(z)$ be a polynomial in 1 complex variable $z$.

Let $\alpha _n$ be randomly chosen from the uniform distribution on the set $\{\alpha \in \mathbf{C}:~|\alpha _n-1|\leq \rho \}$. 

Then the Random Relaxed Newton's method, applied to finding roots of the polynomial $P(z)$, will converge to a root of $P(z)$ for all - except a finite number of exceptions - initial points $z_0$.  

\label{TheoremSumi}\end{thm}

\subsection{Backtracking New Q-Newton's method}
In this section, we recall enough details on bothe theoretical and practical aspects of BNQN.

\subsubsection{Heuristics and some delicate issues} A tendency of Newton's method, applied to find minima of a function $F:\mathbf{R}^m\rightarrow \mathbf{R}$, is that if the initial point $z_0$ is close  to a (non-degenerate) critical point $z^*$ of $F(z)$, then the sequence constructed by Newton's method will converge to $z^*$. Hence, if $z^*$ is a saddle point or local maximum, then Newton's method is not desirable for finding (local) minima. In the special case where $F(z)$ is a quadratic function whose Hessian is invertible, then for every initial point $z_0$, the point $z_1$ constructed by Newton's method will be the unique critical point $z^*=0$.

To overcome this undesirable behaviour of Newton's method, a new variant called New Q-Newton's method (NQN) was recently proposed in \cite{RefTT}. This consists of the following two main ideas: 

- Add a perturbation $\delta ||\nabla F(z)||^{\tau}Id$ to the Hessian $\nabla ^2 F(z)$, where $\tau >0$ is a constant, and $\delta$ is appropriately chosen within a previously chosen set $\{\delta _0,\ldots ,\delta _m\}\subset \mathbf{R}$. Thus, we consider a matrix $A(z)=\nabla ^2F(z)+\delta ||\nabla F(z)||^{\tau}Id$, instead of $\nabla ^2F(z)$ as in the vanilla Newton's method. This has two benefits. First,  it avoids the (minor) difficulty one encounters in Newton's method if $\nabla^2F(z)$ is a singular matrix. Second, more importantly, it turns out that if the $\delta _0,\ldots ,\delta _m$ are randomly chosen, then NQN can avoid saddle points. Since the perturbation $\delta ||\nabla F(z)||^{\tau}Id$ is negligble near non-degenerate local minima, the rate of convergence of NQN is the same as that of Newton's method there.  

- If one mimics Newton's method, in defining $z_{n+1}=z_n-A(z_n)^{-1}.\nabla F(z)$, then the constructed sequence still has the same tendency of convergence to the nearest critical point of $F(z)$. This is remedied in NQN by the following idea: We let $B(z)$ be the matrix with the same eigenvectors as $A(z)$, but whose eigenvalues are all absolute values of the corresponding eigenvalues of $A(z)$. The update rule of NQN is $z_{n+1}=z_n-B(z_n)^{-1}.\nabla F(z_n)$. 

The precise definition of the algorithm NQN is given in the next subsection. If $F(z)$ is a $C^3$ function, then NQN applied to $F(z)$ can avoid saddle points, and has the same rate of convergence near non-degenerate local minima. 

However, NQN has no global convergence guarantee. In \cite{RefT}, a variation of NQN, that is BNQN, was introduced and shown to keep the same good theoretical guarantees as NQN, with the additional bonus that global convergence can be proven for BNQN for very general classes of functions $F(z)$: functions which have at most countably many critical points, or functions which satisfy a Lojasiewicz gradient inequality type. The main idea is to incorporate Armijo's Backtracking line search \cite{RefAr} into NQN. For the readers' convenience, we recall here the idea of Armijo's Backtracking line search. Let $F:\mathbf{R}^m\rightarrow \mathbf{R}$ be a $C^1$ function. Let $z,w\in \mathbf{R}$ such that $<\nabla F(z),w>$ is strictly positive. Then, there exists a positive real number $\gamma$ for which $F(x-\gamma w)-F(x)\leq -<\nabla F(x),\gamma w >/3$. If we choose $\kappa $ by  a backtracking manner (that is, start $\kappa$ from a positive number, and then reduce it exponentially until the above mentioned inequality is satisfied), then the procedure is called Armijo's Backtracking line search.

We notice that there are some delicate issues when doing this incorporation between NQN and Armijo's Backtracking line search, which rely crucially on the fact that we are using a perturbation of the Hessian matrix here. First, in BNQN, the choice of $\delta$ from among $\{\delta _0,\ldots ,\delta _m\}$ is more complicated than NQN. Second, the analog Backtracking line search for Gradient descent is not yet known to be able to avoid saddle points, even though there is a slight variant which can avoid saddle points, and experiments support that the method should be able to avoid saddle point. Third, a priori the learning rate one finds by Armijo's Backtracking line search can be smaller than 1, and hence the rate of convergence can a priori slower than being quadratic.  

The precise definition of BNQN is given in the next subsection. 

\subsubsection{Algorithms: NQN and BNQN} Here we present the basic versions of NQN and BNQN. Many more variations can be found in \cite{RefT}. 

Let $A:\mathbb{R}^m\rightarrow \mathbb{R}^m$ be an invertible {\bf symmetric} square matrix. In particular, it is diagonalisable.  Let $V^{+}$ be the vector space generated by eigenvectors of positive eigenvalues of $A$, and $V^{-}$ the vector space generated by eigenvectors of negative eigenvalues of $A$. Then $pr_{A,+}$ is the orthogonal projection from $\mathbb{R}^m$ to $V^+$, and  $pr_{A,-}$ is the orthogonal projection from $\mathbb{R}^m$ to $V^-$. As usual, $Id$ means the $m\times m$ identity matrix.  

First, we introduce NQN \cite{RefTT}. 

\medskip
{\color{blue}
 \begin{algorithm}[H]
\SetAlgoLined
\KwResult{Find a minimum of $F:\mathbb{R}^m\rightarrow \mathbb{R}$}
Given: $\{\delta_0,\delta_1,\ldots, \delta_{m}\}\subset \mathbb{R}$\  and $\alpha >0$;\\
Initialization: $z_0\in \mathbb{R}^m$\;
 \For{$k=0,1,2\ldots$}{ 
    $j=0$\\
    \If{$\|\nabla f(z_k)\|\neq 0$}{
   \While{$\det(\nabla^2f(z_k)+\delta_j \|\nabla f(z_k)\|^{1+\alpha}Id)=0$}{$j=j+1$}}

$A_k:=\nabla^2f(z_k)+\delta_j \|\nabla f(z_k)\|^{1+\alpha}Id$\\
$v_k:=A_k^{-1}\nabla f(z_k)=pr_{A_k,+}(v_k)+pr_{A_k,-}(v_k)$\\
$w_k:=pr_{A_k,+}(v_k)-pr_{A_k,-}(v_k)$\\
$x_{k+1}:=x_k-w_k$
   }
  \caption{New Q-Newton's method} \label{table:alg}
\end{algorithm}
}
\medskip

BNQN includes a more sophisticated choice of $\delta$ in NQN, together with a combination of Armijo's Backtracking line search. For a symmetric, square real matrix $A$, we define: 
  
  $sp(A)=$ the maximum among $|\lambda |$'s, where $\lambda  $ runs in the set of eigenvalues of $A$, this is usually called the spectral radius in the Linear Algebra literature;
  
  and 
  
  $minsp(A)=$ the minimum among $|\lambda |$'s, where $\lambda  $ runs in the set of eigenvalues of $A$, this number is non-zero precisely when $A$ is invertible.
  
 One can easily check the following more familiar formulas: $sp(A)=\max _{\|e\|=1}\|Ae\|$ and $minsp(A)=\min _{\|e\|=1}\|Ae\|$, using for example the fact that $A$ is diagonalisable.  

We recall that a function $F$ has compact sublevels if for all $C\in \mathbf{R}$ the set $\{z:~F(z)\leq C\}$ is compact. 

\medskip
{\color{blue}
 \begin{algorithm}[H]
\SetAlgoLined
\KwResult{Find a minimum of $F:\mathbb{R}^m\rightarrow \mathbb{R}$}
Given: $\{\delta_0,\delta_1,\ldots, \delta_{m}\} \subset \mathbb{R}$\, $0<\tau $ and $0<\gamma _0\leq 1$;\\
Initialization: $z_0\in \mathbb{R}^m$\;
$\kappa:=\frac{1}{2}\min _{i\not=j}|\delta _i-\delta _j|$;\\
 \For{$k=0,1,2\ldots$}{ 
    $j=0$\\
  \If{$\|\nabla F(z_k)\|\neq 0$}{
   \While{$minsp(\nabla^2F(z_k)+\delta_j \|\nabla F(z_k)\|^{\tau}Id)<\kappa  \|\nabla F(z_k)\|^{\tau}$}{$j=j+1$}}
  
 $A_k:=\nabla^2f(z_k)+\delta_j \|\nabla f(z_k)\|^{\tau}Id$\\
$v_k:=A_k^{-1}\nabla f(z_k)=pr_{A_k,+}(v_k)+pr_{A_k,-}(v_k)$\\
$w_k:=pr_{A_k,+}(v_k)-pr_{A_k,-}(v_k)$\\
$\widehat{w_k}:=w_k/\max \{1,\|w_k\|\}$\\
(If $F$ has compact sublevels, then one can choose $ \widehat{w_k}=w_k$).\\
%(\If{F has compact supports}{ can choose \widehat{w_k}=w_k}.)
$\gamma :=\gamma _0$\\
 \If{$\|\nabla f(z_k)\|\neq 0$}{
   \While{$f(z_k-\gamma \widehat{w_k})-f(z_k)>-\gamma \langle\widehat{w_k},\nabla f(z_k)\rangle/3$}{$\gamma =\gamma /3$}}

$z_{k+1}:=z_k-\gamma \widehat{w_k}$
   }
  \caption{Backtracking New Q-Newton's method} \label{table:alg0}
\end{algorithm}
}
\medskip

Algorithm \ref{table:alg0} has two different versions depending whether the objective function $F$ has compact sublevels or not. Here we introduce a new variant, named BNQN New Variant, with a new parameter, which includes both these versions as special cases. Indeed, $\theta =0$ in BNQN New Variant is the version of BNQN for functions with compact sublevels, and $\theta =1$ in BNQN New Variant is the version of BNQN for general functions. 

\medskip
{\color{blue}
 \begin{algorithm}[H]
\SetAlgoLined
\KwResult{Find a minimum of $F:\mathbb{R}^m\rightarrow \mathbb{R}$}
Given: $\{\delta_0,\delta_1,\ldots, \delta_{m}\} \subset \mathbb{R}$\, $\theta \geq 0$, $0<\tau $ and $0<\gamma _0\leq 1$;\\
Initialization: $z_0\in \mathbb{R}^m$\;
$\kappa:=\frac{1}{2}\min _{i\not=j}|\delta _i-\delta _j|$;\\
 \For{$k=0,1,2\ldots$}{ 
    $j=0$\\
  \If{$\|\nabla F(z_k)\|\neq 0$}{
   \While{$minsp(\nabla^2F(z_k)+\delta_j \|\nabla F(z_k)\|^{\tau}Id)<\kappa  \|\nabla F(z_k)\|^{\tau}$}{$j=j+1$}}
  
 $A_k:=\nabla^2f(z_k)+\delta_j \|\nabla f(z_k)\|^{\tau}Id$\\
$v_k:=A_k^{-1}\nabla f(z_k)=pr_{A_k,+}(v_k)+pr_{A_k,-}(v_k)$\\
$w_k:=pr_{A_k,+}(v_k)-pr_{A_k,-}(v_k)$\\
$\widehat{w_k}:=w_k/\max \{1,\theta \|w_k\|\}$\\
%(\If{F has compact supports}{ can choose \widehat{w_k}=w_k}.)
$\gamma :=\gamma _0$\\
 \If{$\|\nabla f(z_k)\|\neq 0$}{
   \While{$f(z_k-\gamma \widehat{w_k})-f(z_k)>-\gamma \langle\widehat{w_k},\nabla f(z_k)\rangle/3$}{$\gamma =\gamma /3$}}

$z_{k+1}:=z_k-\gamma \widehat{w_k}$
   }
  \caption{Backtracking New Q-Newton's method New Variant} \label{table:alg2}
\end{algorithm}
}
\medskip

\subsubsection{A main theoretical result for finding roots of meromorphic functions by BNQN}

As mentioned previously, BNQN has strong theoretical guarantees for some big classes of functions in any dimension. However, to keep the presentation concise, here we present only one main result relevant to the question pursued in this paper, that of finding roots of meromorphic functions in 1 complex variable. For more general results, the readers can consult \cite{RefT}, \cite{RefTT}. 

\begin{thm} Let $g(z):\mathbf{C}\rightarrow \mathbf{P}^1$ be a non-constant meromorphic function. 
Define a function $F:\mathbf{R}^2\rightarrow [0,+\infty]$ by the formula $F(x,y)=|g(x+iy)|^2/2$. 

Given an initial point $z_0\in \mathbf{C}$, which is not a pole of $g$, we let $\{z_n\}$ be the sequence constructed by BNQN New Variant applied to the function $F$ with initial point $z_0$. If the objective function has compact sublevels (i.e. for all $C\in \mathbf{R}$ the set $\{(x,y)\in \mathbf{R}^2:~F(x,y)\leq C\}$ is compact), we choose $\theta \geq 0$, while in general we choose $\theta >0$. 

1) Any critical point of $F$ is a root of $g(z)g'(z)=0$. 

2) If $z^*$ is a cluster point of $\{z_n\}$ (that is, if it is the limit of a subsequence of $\{z_n\}$), then $z^*$ is a critical point of $F$. 

3) If $F$ has compact sublevels, then $\{z_n\}$ converges. 

4) Assume that the parameters $\delta _0,\delta _1,\delta _2$ in BNQN New Variant are randomly chosen. Assume also that $g(z)$ is generic, in the sense that $\{z\in \mathbf{C}:~g(z)g"(z)=g'(z)=0\}=\emptyset$. There exists an exceptional set $\mathcal{E}\subset \mathbf{C}$ of zero Lebesgue measure so that if $z_0\in \mathbf{C}\backslash \mathcal{E}$, then $\{z_n\}$ must satisfy one of the following two options: 

Option 1: $\{z_n\}$ converges to a root $z^*$ of $g(z)$, and if $\gamma _0=1$ in the algorithm then the rate of convergence is quadratic. 

Option 2: $\lim _{n\rightarrow\infty}|z_n| =+\infty$. 

\label{TheoremMeromorphic}\end{thm}

For BNQN, Theorem \ref{TheoremMeromorphic} is proven in \cite{RefTT}. It is easy to modify the proof for BNQN New Variant. An example for which part 3 of Theorem \ref{TheoremMeromorphic} can apply is when $g$ is a polynomial or $g=P/Q$ where $P,Q$ are polynomials and $P$ has bigger degree than $Q$ (a special case is when $Q=P'$, in which case the zeros of $g$ are exactly that of $P$, with the advantage that they all have multiplicity $1$).  

\subsubsection{Implementation details} An implementation in Python of BNQN accompanies the paper \cite{RefT}. The implementation is flexible in that one does not need precise values of the gradient and Hessian of the function $F$, but approximate values are good enough. Experiments also show that the performance of BNQN is very stable, with respect to its parameters and to the values of the objective function and its gradient and Hessian matrix.  

\section{Proofs of Theorems \ref{TheoremMain1} and \ref{TheoremMain2}}
In this section, we first prove Theorem \ref{TheoremMain2}, and then use it to prove Theorem \ref{TheoremMain1}. 

\begin{proof}[Proof of Theorem \ref{TheoremMain2}]
    By induction, it suffices to show that under the assumptions of the theorem, we have $z_1'=A^{-1}z_1$.

    As in the proof of linear conjugacy invariance for Newton's method, we have the following formulas relating the gradients and the Hessians of $F$ and $G$:
	
\begin{itemize}
      \item $\nabla G(z)=A^T \nabla F(Az)=cR^{-1}\nabla F(Az)$
      \item $\nabla^2G(z)=A^T \nabla^2 F(Az) A=c^2R^{-1}\nabla^2 F(Az) R$. 
\end{itemize}
From this, we have that $||\nabla G(z)||=c||\nabla F(Az)||$. Moreover, $(\lambda ,v)\in \mathbf{R}\times \mathbf{R}^m$ is an eigenvalue-eigenvector pair of $\nabla ^2F(Az)$, then $(c^2\lambda ,R^{-1}v)$ is an eigenvalue-eigenvector pair of $\nabla ^2G(z)$.  

If $\kappa =\min _{i\not=j}|\delta _i-\delta _j|/2$ and $\kappa '=\min _{i\not=j}|\delta '_i-\delta '_j|/2$, then we have the relation: 
\begin{eqnarray*}
\kappa '=c^{2-\tau }\kappa . 
\end{eqnarray*}

Therefore, we have for every $\delta =\delta _j$ and $\delta '=\delta _j'$: 
\begin{eqnarray*}
&&\nabla ^2G(z)+\delta _j'||\nabla G(z)||^{\tau}Id\\
&=&c^2R^{-1}.\nabla ^2F(Az).R+c^2\delta _j||\nabla F(Az)||^{\tau}Id\\
&=&c^2R^{-1}(\nabla ^2F(Az)+\delta _j||\nabla F(Az)||^{\tau}Id).R 
\end{eqnarray*}

It follows that $\delta _j$ satisfies $minsp (\nabla ^2F(Az)+\delta _j||\nabla F(Az)||^{\tau}Id)\geq \kappa ||\nabla F(Az)||^{\tau}$ if and only if $\delta _j'$ satisfies  $minsp (\nabla ^2G(z)+\delta '_j||\nabla G(z)||^{\tau}Id)\geq \kappa '||\nabla G(z)||^{\tau}$. 

Hence, if $A_0=\nabla ^2F(Az)+\delta _j||\nabla F(Az)||^{\tau}Id$ is the matrix in the BNQN New Variant algorithm for $F$ (at the point $Az$), then $A_0'=\nabla ^2G(z)+\delta '_j||\nabla G(z)||^{\tau}Id$ is the matrix in the BNQN New Variant algorithm for $G$ (at the point $z$). 

Next, we need to relate the two vectors $v_0$ (for F) and $v_0'$ (for G), in the BNQN New Variant algorithm. We have 
\begin{eqnarray*}
v_0'&=&(A_0')^{-1}.\nabla G(z)\\
&=&(c^2R^{-1}(\nabla ^2F(Az)+\delta _j||\nabla F(Az)||^{\tau}Id).R)^{-1}.(cR^{-1}.\nabla F(Az))\\
&=&c^{-1}R^{-1}.(\nabla ^2F(Az)+\delta _j||\nabla F(Az)||^{\tau}Id).\nabla F(Az)\\
&=&A^{-1}.v_0. 
\end{eqnarray*}

From this, it follows easily, from the relation between $A_0$ and $A_0'$, that 
\begin{eqnarray*}
w_0'&=&pr_{A',+}v_0'-pr_{A',-}v_0'\\
&=&A^{-1}.w_0. 
\end{eqnarray*}

We note that for general $c$: 
\begin{eqnarray*}
<w_0',\nabla G(z)>=<c^{-1}R^{-1}.w_0,cR^{-1}\nabla F(Az) >=<w_0,\nabla F(Az)>.  
\end{eqnarray*}

Also, with $\theta '=c \theta$, we have $\max \{1,\theta '||w_0'||\}=\max \{1, \theta ||w_0||\}$. 

Therefore, if $\widehat{w_0}=w_0/\max \{1,\theta ||w_0||\}$ and $\widehat{w_0'}=w_0'/\max \{1,\theta '||w_0'||\}$, we have 
\begin{eqnarray*}
<\widehat{w_0'},\nabla G(z)>=<\widehat{w_0},\nabla F(Az)>.  
\end{eqnarray*}

Therefore, the learning rates $\gamma$ and $\gamma'$, found by applying BNQN New Variant to $F$ (at the point $z_0=Az_0'$) and $G$ (at the point $z_0'$) are the same. 

This implies that 
\begin{eqnarray*}
z_1'=z_0'-\gamma '\widehat{w_0'}=A^{-1}.z_0-\gamma A^{-1}.\widehat{w_0}=A^{-1}.(z_0-\gamma \widehat{w_0})=A^{-1}.z_1,    
\end{eqnarray*}
as wanted. 
\end{proof}

We next present a simple example showing that if the linear conjugacy $A$ is not of the form $cR$, where $c\in \mathbf{R}$ and $R$ is a unitary operator, then the constructed sequences $\{z_n\}$ and $\{z_n'\}$ may not satisfy the relation $z_n'=A^{-1}z_n$ for all $n$. 

\begin{example} Let $F(x,y)=xy$, and $A$ be the matrix
\[\begin{bmatrix}
1&1\\
0&1
\end{bmatrix}.\]

$A$ is not of the form $cR$ where $c\in \mathbf{R}$ and $R^TR=Id$, since $A^T$ is the matrix 
matrix
\[\begin{bmatrix}
1&0\\
1&1
\end{bmatrix}\]
and $AA^T$ is the matrix
\[\begin{bmatrix}
2&1\\
1&1
\end{bmatrix}.\]
We can also compute that $A^{-1}$ is the matrix
\[\begin{bmatrix}
1&-1\\
0&1
\end{bmatrix}.\]

Define $G(x,y)=F(A.(x,y))=(x+y)y$.

We can compute $\nabla F(x,y)=(y,x)$, and $\nabla G(x,y)=(y, x+2y)$. 

The Hessian $\nabla ^2F(x,y)$ is 
\[\begin{bmatrix}
0&1\\
1&0
\end{bmatrix}, \]
and its inverse is 
\[\begin{bmatrix}
0&1\\
1&0
\end{bmatrix}. \]

The Hessian $\nabla ^2G(x,y)$ is 
\[\begin{bmatrix}
0&1\\
1&2
\end{bmatrix},  \]
and its inverse is 
\[\begin{bmatrix}
-2&1\\
1&0
\end{bmatrix}. \]

$\nabla ^2F$ has two eigenvalue-eigenvector pairs $(\lambda _1,u_1)$ and $(\lambda _2,u_2)$ where $\lambda _1=-1$, $u_1=(-1,1)$ and $\lambda _2=1$, $u_2=(1,1)$. 

$\nabla ^2G$ has two eigenvalue-eigenvector pairs $(\lambda '_1,u'_1)$ and $(\lambda '_2,u'_2)$ where $\lambda '_1=1+\sqrt{2}$, $u'_1=(-1+\sqrt{2},1)$ and $\lambda '_2=1-\sqrt{2}$, $u'_2=(-1-\sqrt{2},1)$. 

We next compute the terms in the definition of BNQN New Variant (where for simplicity we choose $\delta =\delta '=0$, and hence $A_0=\nabla ^2F$ and $A_0'=\nabla ^2G$, but the argument can be easily modified to fit other values of $\delta$):
\begin{eqnarray*}
v_0&=&(\nabla ^2F(A.(x,y)))^{-1}.\nabla F(A.(x,y))=(x+y,y)\\
&=&-x.(-1,1)/2 + (x+2y).(1,1)/2\\ 
v_0'&=&(\nabla ^2G(x,y))^{-1}.\nabla G(x,y)=(x,y).
\end{eqnarray*}

Therefore, 
\begin{eqnarray*}
w_0&=&x.(-1,1)/2 +(x+2y)(1,1)/2 =(y,x+y),\\
A^{-1}.w_0&=&(-x,x+y).  
\end{eqnarray*}

We show that $w_0'$ cannot be a multiple of $A^{-1}.w_0$. Indeed, if $w_0'=a.(A^{-1}.w_0)$, then by the definition of $w_0'$, we must have 
\begin{eqnarray*}
x(-1+\sqrt{2})+y
&=&<v_0',u_1'>\\
&=&a<A^{-1}.w_0,u_1'>\\&=&a((-x)(-1+\sqrt{2})+x+y),\\
x(-1-\sqrt{2})+y&=&<v_0',u_2'>\\
&=&-a<A^{-1}.w_0,u_2'>\\
&=&-a((-x)(-1-\sqrt{2})+x+y). 
\end{eqnarray*}
Thus, dividing side by side the above two equations, assuming all the terms are non-zero, we obtain a relation between $x,y$, if we assume that $w_0'=aA^{-1}.w_0$: 
\begin{eqnarray*}
    \frac{x(-1+\sqrt{2})+y}{x(-1-\sqrt{2})+y}=-\frac{(-x)(-1+\sqrt{2})+x+y}{(-x)(-1-\sqrt{2})+x+y}. 
\end{eqnarray*}
For a generic choice of $(x,y)$, the above equality cannot be satisfied. Therefore, for a generic choice of $(x,y)$, $w_0'\not= a.A^{-1}.w_0$ for any real number $a$. 

Since $z_1=z_0-\gamma w_0$ and $z_1'=z_0'-\gamma 'w_0'$, and $z_0'=A^{-1}.z_0$, the above implies that $z_1'\not= A^{-1}.z_1$. 

\end{example}

We are now ready to prove Theorem \ref{TheoremMain1}.

\begin{proof}[Proof of Theorem \ref{TheoremMain1}] By Theorem \ref{TheoremMain2}, it suffices to consider two special cases: $f(z)=z^2$ and $f(z)=z^2-1$. 

In this case, $F(x,y)=|f(x+iy)|^2/2$ has compact sublevels. Hence, Theorem \ref{TheoremMeromorphic} implies that for every initial point $z_0$, the sequence $\{z_n\}$ constructed by BNQN New Variant will converge to a point $z^*$, which must be a root of $f(z)f'(z)$. 

1) The case $f(z)=z^2$: In this case, $f(z)f'(z)=2z^3$ has only one root (with multiplicity 3). Hence, $z^*=0$, as wanted. 

2) The case $f(z)=z^2-1$: $f$ has 2 roots $z_1^*=-1$ and $z_2^*=1$, and has a critical point $z_3^*=0$. Therefore, $z^*$ must be one of these 3 points. 

The line $L$ is the $y$-axis: $L=\{(x,y)\in \mathbf{R}^2:~x=0\}$.  

a) We first show that if $z_0\in H_1=\{(x,y)\in \mathbf{R}^2:~x>0\}$, then the limit point $z^*$ of the constructed sequence $\{z_n\}$ is $1$. (At the end of the proof, we will explain how Theorem \ref{TheoremMain2} can be used to deal with the case $z_0\in H_2=\{(x,y)\in \mathbf{R}^2:~x<0\}$.) 

To this end, we will prove the following two claims: 

{\bf Claim A:} If $z_n=(x_n,y_n)$, then $x_n>0$ for all $n$. 

{\bf Claim B:} (Stable manifold for the saddle point $0$) If $z_n\rightarrow (0,0)$ then $x_n=0$ for $n$ large enough. 

We now explain how the above two Claims shows that $z^*=1$. Indeed, since $z^*$ is either $-1$, $0$ or $1$, and by Claim A we have $x_n>0$ for all $n$, we have $z^*$ cannot be $-1$, and must be either $1$ or $0$. Then, Claim B implies that $z^*$ cannot be $0$. Hence $z^*=1$. 

We now proceed to proving Claim A. It suffices to show that if $z_0=(x,y)$ has $x>0$, then $z_1=(x_1,y_1)$ has $x_1>0$. 

We have
$$
F(x, y)=\frac{1}{2}||f(x+iy)||^2=\frac{1}{2}[(x^2-y^2-1)^2+4x^2y^2];
$$

$$
 \triangledown F=
 \left[ \begin{array}{ccc}
 	2(x^2+y^2-1)x \\
 	2(x^2+y^2+1)y \\
 \end{array} \right];
$$
$$ 
 \triangledown ^2 F=
\left(
\begin{array}{cc}
	6x^2+2y^2-2 & 4xy \\
	4xy& 2x^2+6y^2+2
\end{array}
\right).
$$

Since the function $F(x,y)$ is invariant under the rotation $(x,y)\mapsto (x,-y)$, by Theorem \ref{TheoremMain2}  we can assume that $z_0=(x,y)$ where $x>0, y\geq 0$. We treat separately the two cases: $y>0$ and $y=0$.  

\underline{\bf Case 1:} $x,y>0$. We will show that for $z_1=(x_1,y_1)$ obtained from $z_0=(x,y)$ by BNQN New Variant, then $x_1>0$. Assume otherwise that $x_1\leq 0$, we will arrive at a contradiction. 

The eigenvalues of $ \triangledown ^2 F$ are
$$\lambda_1=2(2x^2-\sqrt{x^4+2x^2y^2-2x^2+y^4+2y^2+1}+2y^2)$$ 
and 
$$\lambda_2=2(2x^2+\sqrt{x^4+2x^2y^2-2x^2+y^4+2y^2+1}+2y^2).$$ 
The eigenvector of $ \triangledown ^2 F$ for $\lambda_1$ is
$$
u_1=
\left(
\begin{array}{cc}
	x^2-y^2-1-\sqrt{1-2x^2+x^4+2y^2+2x^2y^2+y^4} \\
	2xy
\end{array}
\right)
.$$

Note that the term under the square root can also be written
as $$
S:=\sqrt{(1-x^2+y^2)^2+4x^2y^2}.
$$

Then we can write:

$$
u_1=
\left(
\begin{array}{cc}
	x^2-y^2-1-S \\
	2xy
\end{array}
\right)
$$ 

The eigenvector of $ \triangledown ^2 F$ for $\lambda_2$ is
\begin{eqnarray*}
u_2&=&
\left(
\begin{array}{cc}
	x^2-y^2-1+\sqrt{1-2x^2+x^4+2y^2+2x^2y^2+y^4} \\
	2xy
\end{array}
\right)
\\&=&
\left(
\begin{array}{cc}
	x^2-y^2-1+S \\
	2xy
\end{array}
\right).
\end{eqnarray*}

{\bf Claim 1:} If $x>0$ and $y>0$, then there are positive real numbers $a_1$ and $a_2$ so that $\nabla F=a_1u_1+a_2u_2$. 

\begin{proof}[Proof of Claim 1] 
It suffices to check that $<\nabla F,u_1>$ and $<\nabla F,u_2>$ are positive numbers. 

We have 
\begin{equation*}
\begin{aligned}
    <u_1, \nabla F(z_0)>=&[-(1-x^2+y^2)-\sqrt{1-2x^2+x^4+2y^2+2x^2y^2+y^4}]\cdot\\
    & 2(x^2+y^2-1)x
    +2xy\cdot 2(x^2+y^2+1)y\\
    =&2x[(x^2-1)^2+2y^2+2x^2y^4+y^4]-\\
    &2x(x^2+y^2-1)\sqrt{1-2x^2+x^4+2y^2+2x^2y^2+y^4}\\
    =&2x\sqrt{1-2x^2+x^4+2y^2+2x^2y^2+y^4}\cdot\\
    &[\sqrt{1-2x^2+x^4+2y^2+2x^2y^2+y^4}-(x^2+y^2-1)]
\end{aligned}
\end{equation*}
We know 
\begin{equation*}
\begin{aligned}
&1-2x^2+x^4+2y^2+2x^2y^2+y^4-(x^2+y^2-1)^2\\
=&4y^2>0.
\end{aligned}
\end{equation*}
Hence, $$<u_1, \nabla F(z_0)> ~ >0.$$
\begin{equation*}
\begin{aligned}
    <u_2, \nabla F(z_0)>=&[-(1-x^2+y^2)+\sqrt{1-2x^2+x^4+2y^2+2x^2y^2+y^4}]\cdot\\
    & 2(x^2+y^2-1)x
    +2xy\cdot 2(x^2+y^2+1)y\\
   =&2x\sqrt{1-2x^2+x^4+2y^2+2x^2y^2+y^4}\cdot\\
    &[\sqrt{1-2x^2+x^4+2y^2+2x^2y^2+y^4}+(x^2+y^2-1)]
\end{aligned}
\end{equation*}

Thus, $$<u_2, \nabla F(z_0)> ~ >0.$$

\end{proof}

Let, in BNQN New Variant, $A_0=\nabla ^2F(x,y)+\delta ||\nabla F(x,y)||^{\tau}$, where $A_0$ is invertible and has two eigenvalues $\lambda _1'$, $\lambda _2'$.  By definition $z_1=z_0-\gamma w_0$, for some $\gamma >0$, and $w_0$ is given by the formula: 
\begin{eqnarray*}
w_0=\frac{a_1}{|\lambda '_1|}u_1+\frac{a_2}{|\lambda _2'|}u_2. 
\end{eqnarray*}
Therefore, for $U_1=z_0-z_1=\gamma w_0$, by Claim 1 we can write $U_1=b_1u_1+b_2u_2$ for $b_1,b_2>0$. 

Let $\hat{z}=(0,\hat{y})$ be the intersection point between the line connecting the two points $z_0$ and $z_1$, and the line $L=\{(x,y):~x=0\}$. Since $x>0$ and $x_1\leq 0$, the point $\hat{z}$ is in between $z_0$ and $z_1$.  

{\bf Claim 2:} If $x_1\leq 0$, then \begin{equation}\label{eq8}
y_1\leq \hat{y}\leq \frac{x^2+y^2-1-\sqrt{(x^2+y^2-1)^2+4y^2}}{2y}.
\end{equation}
\begin{proof}[Proof of Claim 2]
Let $\hat{U}=z_0-\hat{z}$. Then $\hat{U}$ has the same direction as $U_1$. Therefore, 

\begin{equation*}
\begin{aligned}
&<\hat{U}, u_1>\\
=&x \cdot [-(1-x^2+y^2)-\sqrt{1-2x^2+x^4+2y^2+2x^2y^2+y^4}]+2xy(y-\hat{y})>0\\
\end{aligned}
\end{equation*}  
Thus, since $x>0$, we have
\begin{equation*}
\begin{aligned}
\hat{y}<&\frac{x^2+y^2-1-\sqrt{1-2x^2+x^4+2y^2+2x^2y^2+y^4}}{2y}\\
=&\frac{x^2+y^2-1-\sqrt{(x^2+y^2-1)^2+4y^2}}{2y}.\\
%<&0.
\end{aligned}
\end{equation*} 

Similarly, 
\begin{equation*}
\begin{aligned}
&<\hat{U}, u_2>\\
=&x \cdot [-(1-x^2+y^2)+\sqrt{1-2x^2+x^4+2y^2+2x^2y^2+y^4}]+2xy(y-\hat{y})>0\\
\end{aligned}
\end{equation*}  

Hence,  
$$\hat{y}\leq\frac{x^2+y^2-1+\sqrt{1-2x^2+x^4+2y^2+2x^2y^2+y^4}}{2y}.$$

Therefore, $\hat{y}$ is smaller than the minimum of the two expressions, which is 
$$\frac{x^2+y^2-1-\sqrt{(x^2+y^2-1)^2+4y^2}}{2y}.$$

Since $x>0,y>0$ and $\hat{y}\leq 0$, it follows that $y_1\leq \hat{y}$. 
\end{proof}

Now we explore the RHS term in the upper bound for $\hat{y}$ in Claim 2. 

{\bf Claim 3:} If $x,y>0$ and $x^2+y^2-1\leq t$, then 
\begin{eqnarray*}
\frac{x^2+y^2-1-\sqrt{(x^2+y^2-1)^2+4y^2}}{2y}\leq \frac{t-\sqrt{t^2+4y^2}}{2y}. 
\end{eqnarray*}

In particular: 

a) If  $0<x\leq 1$, $0< y\leq 1$, then 
\begin{eqnarray*}
\frac{x^2+y^2-1-\sqrt{(x^2+y^2-1)^2+4y^2}}{2y}\leq \frac{1-\sqrt{5}}{2}y\leq -0.61y. 
\end{eqnarray*}

b) If $x^2+y^2\leq 1$, then 
\begin{eqnarray*}
\frac{x^2+y^2-1-\sqrt{(x^2+y^2-1)^2+4y^2}}{2y}\leq -1. 
\end{eqnarray*}

\begin{proof}[Proof of Claim 3] We define $s=x^2+y^2-1\in \mathbf{R}$, and $a=4y^2>0$. Then 
\begin{eqnarray*}
\frac{x^2+y^2-1-\sqrt{(x^2+y^2-1)^2+4y^2}}{2y}= \frac{s-\sqrt{s^2+a}}{2y}. 
\end{eqnarray*}
We consider the function $h(s)=s-\sqrt{s^2+a}$. Its derivative is $h'(s)=1-s/\sqrt{s^2+a}$ which is always positive, due to $a>0$. Therefore,  $h(s)$ is an increasing function. Since $y>0$, this implies the desired inequality 
\begin{eqnarray*}
\frac{x^2+y^2-1-\sqrt{(x^2+y^2-1)^2+4y^2}}{2y}\leq \frac{t-\sqrt{t^2+4y^2}}{2y}. 
\end{eqnarray*}

a) Since $0<x\leq 1$, we have $x^2+y^2-1\leq t=y^2$. By the above inequality we get
\begin{eqnarray*}
\frac{x^2+y^2-1-\sqrt{(x^2+y^2-1)^2+4y^2}}{2y}\leq \frac{y^2-\sqrt{y^4+4y^2}}{2y}. 
\end{eqnarray*}
Since $0<y\leq 1$, we have $y^4+4y^2\geq 5y^4$. Therefore
\begin{eqnarray*}
\frac{y^2-\sqrt{y^4+4y^2}}{2y}\leq \frac{y^2-\sqrt{5y^4}}{2y}=\frac{1-\sqrt{5}}{2}y, 
\end{eqnarray*}
and we obtain the desired inequality. 

b) Since $x^2+y^2\leq 1$, we have $x^2+y^2-1\leq t=0$. By the above inequality, we get
\begin{eqnarray*}
\frac{x^2+y^2-1-\sqrt{(x^2+y^2-1)^2+4y^2}}{2y}\leq \frac{0-\sqrt{0+4y^2}}{2y}=-1,  
\end{eqnarray*}
as wanted. 
\end{proof}

We next obtain a useful lower bound for $-2y\hat{y}$, which is needed later. 

{\bf Claim 4:} Under the same assumptions (that is $x,y>0$ and $x_1\leq 0$) we have: 
\begin{eqnarray*}
-2yy_1\geq -2y\hat{y}\geq 2y+1-x^2-y^2. 
\end{eqnarray*}
\begin{proof}[Proof of Claim 4]

By Claim 2, we have 
\begin{eqnarray*}
-2yy_1\geq -2y\hat{y}&\geq&\sqrt{(x^2+y^2-1)^2+4y^2}+1-x^2-y^2\\
&\geq& \sqrt{4y^2}+1-x^2-y^2\\
&=&2y+1-x^2-y^2. 
\end{eqnarray*}
\end{proof}

Again, we recall the assumption that $z_0=(x,y)$ with $x,y>0$, and $z_1=(x_1,y_1)$ with $x_1\leq 0$ is obtained from $z_0$ by BNQN New Variant. Recall  $U_1=z_0-z_1$ and $\hat{U}=z_0-\hat{z}$, where $\hat{z}=(0,\hat{y})$ is the intersection point between the line connecting $z_0$ and $z_1$, with the line $L$. We will now  show that 
\begin{equation}
F(z_1)-F(z_0)+<U_1,\nabla F(z_0)>/3 >0, 
\label{EquationContradictionArmijoInequality}\end{equation}
and obtain a contradiction to the requirement on Armijo's Backtracking line search in BNQN New Variant. 

We define also $U=\hat{z}-z_1$. We can write 
\begin{eqnarray*}
F(z_1)-F(z_0)+<U_1,\nabla F(z_0)>/3=I+II+III,
\end{eqnarray*}
where 
\begin{itemize}
	\item $I= F(0,\hat{y})-F(x,y)+\frac{1}{3}\bigg<
	\left( \begin{array}{ccc}
		x \\
		y-\hat{y} \\
	\end{array} \right),\nabla F(x,y)\bigg>.$\\
	\item $II=F(x_1,y_1)-F(0,\hat{y}).$\\
	\item $III=\frac{1}{3}\bigg<\left( \begin{array}{ccc}
		-x_1 \\
		\hat{y}-y_1 \\
	\end{array} \right),\nabla F(x,y)\bigg>.$\\
\end{itemize}

We have the following explicit formulas for I, II and III: 
\begin{itemize}
	\item $I= \frac{1}{2}\hat{y}^4+\hat{y}^2+\frac{1}{6}x^4+\frac{1}{6}y^4+\frac{1}{3}x^2y^2+\frac{1}{3}x^2-\frac{2}{3}yx^2\hat{y}-\frac{2}{3}y^3\hat{y}-\frac{1}{3}y(y+2\hat{y}).$\\
	\item $II= \frac{1}{2}(y_1^4-\hat{y}^4)+(y_1^2-\hat{y}^2)+\frac{1}{2}x_1^4+x_1^2y_1^2-x_1^2.$\\
\item $III= \frac{1}{3}<U,\nabla F(x,y) >$.
\end{itemize}

{\bf Claim 5:} We have $I>0$ and $III\geq 0$. 
\begin{proof}[Proof of Claim 5] 
Since $U$ is either 0 or has the same direction as $U_1$, we have as  in the proof of Claim 2, that $<U, \nabla F(x,y)> ~\geq 0$, which implies that $III\geq 0$. 

Now we show that $I>0$. From the expression for $I$ and Claim 4, we have  
\begin{eqnarray*}
I&>&\frac{x^4}{6}+\frac{x^2}{3}+\frac{y^4}{6}-\frac{y^2}{3}-\frac{2y\hat{y}}{3}\\
&\geq&\frac{x^4}{6}+\frac{x^2}{3}+\frac{y^4}{6}-\frac{y^2}{3}+\frac{2y+1-x^2-y^2}{3}\\
&=&\frac{1}{6}(x^4+y^4+4y+2-4y^2). 
\end{eqnarray*}
Since $y^4+4y^2-4y^2\geq 0$ for all $y\geq 0$ (see Claim 6 below), we have $I>0$ as wanted. 

\end{proof}

\underline{Subcase 1.1: $x\geq 1$ or $y\geq 1$}: We will show that $I+II+III>0$ and hence obtain the wanted contradiction (\ref{EquationContradictionArmijoInequality}). 

Since $III\geq 0$ by Claim 5, it is enough to show that $I+II> 0$. We have 
\begin{eqnarray*}
 I+II> \frac{x^4}{6}+\frac{x^2}{3}+\frac{y^4}{6}-\frac{y^2}{3}-\frac{2\hat{y}y}{3}+\frac{1}{2}x_1^4-x_1^2,
\end{eqnarray*}
and hence by Claim 4: 
\begin{eqnarray*}
 I+II&>& \frac{x^4}{6}+\frac{x^2}{3}+\frac{y^4}{6}-\frac{y^2}{3}+\frac{2y+1-x^2-y^2}{3}+\frac{1}{2}x_1^4-x_1^2\\
 &=&\frac{1}{6}(x^4+y^4+4y+2-4y^2)+\frac{1}{2}x_1^4-x_1^2. 
\end{eqnarray*}
Hence, $I+II> 0$ provided we can show that
\begin{eqnarray*}
\frac{1}{6}(x^4+y^4+4y+2-4y^2)\geq 1/2, 
\end{eqnarray*}
when $x\geq 1$ or $y\geq 1$. This relies on the following inequalities: 

{\bf Claim 6:} 

a) For all $y\geq 0$, we have $y^4+4y-4y^2\geq 0$. 

b) For all $y\geq 1$, we have $y^4+4y-4y^2\geq 1$. 
\begin{proof}[Proof of Claim 6]

One can easily check this by finding the critical points of the function $h(y)=y^4+4y-4y^2$ and evaluate the function value at the critical points and at $0$ and $1$. Indeed, the derivative of the function is $h'(y)= 4(y^3-2y+1)$, and it has 3 roots $1$, $-(1+\sqrt{5})/2$ and $(\sqrt{5}-1)/2$. We have $-(1+\sqrt{5})/2<0$ and hence needs not be considered. The point $(\sqrt{5}-1)/2$ is in the interval $[0,1]$, and $h((\sqrt{5}-1)/2)$ is about $1.09$. Also, $h(0)=0$ and $h(1)=1$. 
\end{proof}

By Claim 6, if $x\geq 1$ or $y\geq 1$, we have $x^4+y^4+4y-4y^2+2\geq 1+2=3$, and hence 
\begin{eqnarray*}
\frac{1}{6}(x^4+y^4+4y+2-4y^2)\geq 1/2, 
\end{eqnarray*}
as wanted. 

\underline{Subcase 1.2: $0<x,y\leq 1$ and  $x^2+y^2\geq 1$:} Again, we will show that $I+II> 0$, and hence the inequality (\ref{EquationContradictionArmijoInequality}) holds. 

We know
$$I+II>\hat{y}^2-\frac{y^2}{3}-\frac{2}{3}y\hat{y}+\frac{x^2}{3}+\frac{x^4+2x^2y^2+y^4}{6} +\frac{x_1^4}{2} -x_1^2.$$

By part a of Claim 3, $\hat{y}\leq \frac{1-\sqrt{5}}{2}y<0$. Hence, $\hat{y}^2> \frac{(1-\sqrt{5})^2}{4}y^2$. Thus, we have $\hat{y}^2-\frac{y^2}{3}-\frac{2}{3}y\hat{y}>\big[\frac{(1-\sqrt{5})^2}{4}-\frac{1}{3}+\frac{2}{3}\frac{\sqrt{5}-1}{2}\big]y^2=\frac{5-\sqrt{5}}{6}y^2>\frac{1}{3}y^2.$ By using $x^2+y^2\geq 1$,  it follows that
$$I+II>\frac{1}{3}+\frac{1}{6}+\frac{1}{2}x_1^4-x_1^2=\frac{1}{2}+\frac{1}{2}x_1^4-x_1^2\geq 0,$$
as wanted. 

\underline{Subcase 1.3: $0<x,y\leq 1$ and $x^2+y^2\leq 1$:} We will show that $II\geq 0$, and hence the inequality (\ref{EquationContradictionArmijoInequality}) holds by Claim 5. 

We have $II\geq x_1^2y_1^2-x_1^2$. By Claim 2 and part b of Claim 3, we have $y_1\leq \hat{y}\leq -1$. Hence $y_1^2\geq 1$, ande $x_1^2y_1^2-x_1^2\geq 0$, as wanted. 

To finish the proof of Claim A, we explore the remaining case: 

\underline{\bf Case 2: $x>0,y=0$:} Again, we will show that if $z_1=(x_1,y_1)$ is obtained from $z_0=(x,0)$ by BNQN New Variant, where $x>0$, then $x_1>0$. We can assume that $x\not= 1$, since $1$ is  a critical point of $F$, and hence when $z_0=(1,0)$ we also have $z_n=(1,0)$ for all $n$.  

We can explicitly calculate 
$$
 \triangledown F=
 \left[ \begin{array}{ccc}
 	2(x^2-1)x \\
 	0 \\
 \end{array} \right];
$$
$$ 
 \triangledown ^2 F=
\left(
\begin{array}{cc}
	6x^2-2 & 0 \\
	0& 2x^2+2
\end{array}
\right).
$$

We note that $\nabla F(z_0)$ itself is an eigenvector of $\nabla ^2F(z_0)$, with the eigenvalue $2(x^2-3)$. Therefore, for $A_0=\nabla ^2F(z_0)+\delta ||\nabla F(z_0)||^{\tau}Id$, which has one non-zero eigenvalue $\lambda _1'$ with eigenvector $\nabla F(z_0)$, we have $w_0=a\nabla F(z_0)$ for some $a>0$. 

\underline{Subcase 2.1: $0<x<1$, $y=0$:} In this case, by definition of BNQN New Variant, we have $z_1=z_0-\gamma w_0=z_0-a \gamma \nabla F(z_0)$, where $\gamma >0$, $a>0$ and the $x$-coordinate of $\nabla F(z_0)$ is $<0$. It follows that the $x$-coordinate of $z_1$ is bigger than the $x$-coordinate of $z_0$, and hence $x_1>0$, as wanted. 

\underline{Subcase 2.2: $x>1$, $y=0$:} In this case, we will show, similar to Case 1, that if $x_1\leq 0$ then the inequality (\ref{EquationContradictionArmijoInequality}) holds, and hence is a contradiction to the definition of BNQN New Variant. 

We always have $F(z_1)\geq 0$, hence to prove (\ref{EquationContradictionArmijoInequality}) it suffices to show that $F(z_0)-<U,\nabla F(z_0)>/3 <0$, for $U=(x-x_1,0)$. 

Explicit calculation shows that 
\begin{eqnarray*}
F(z_0)-\frac{1}{3}<U,\nabla F(z_0)>&=&\frac{1}{2}(x^2-1)^2-\frac{1}{3}(x-x_1)2(x^2-x)x\\
&=&\frac{1}{2}(x^2-1)[(x^2-1)-\frac{4}{3}(x-x_1)x]\\
&=&\frac{1}{2}(x^2-1)(-1-\frac{1}{3}x^2+\frac{4}{3}xx_1).
\end{eqnarray*}
The last term 
\begin{eqnarray*}
\frac{1}{2}(x^2-1)(-1-\frac{1}{3}x^2+\frac{4}{3}xx_1)
\end{eqnarray*}
is $<0$, whenever $x>1$ and $x_1\leq 0$. Hence 
\begin{eqnarray*}
F(z_0)-\frac{1}{3}<U,\nabla F(z_0)><0,  
\end{eqnarray*}
as wanted. 

Therefore, Claim A is completely proven. We now proceed to proving Claim B. 

To this end, it suffices to show that if $\epsilon >0$ is small enough, $z_0=(x,y)$ is in the domain $\{(x,y):~0<x<\epsilon, 0\leq y<\epsilon \}$, and $z_1=(x_1,y_1)$ is constructed from $z_0$ by BNQN New Invariant, then $x_1>x_0$. 

We note that this has been proven for the case $0<x<1, y=0$ in the proof of Claim A for Subcase 2.1. Therefore, we only need to consider $z_0\in B_{\epsilon}=\{(x,y):~0<x<\epsilon, 0< y<\epsilon \}$. We have: 

$$
 \triangledown F=
 \left[ \begin{array}{ccc}
 	2(x^2+y^2-1)x \\
 	2(x^2+y^2+1)y \\
 \end{array} \right] \approx
 \left(
\begin{array}{cc}
	-2x \\
	2y
\end{array}
\right);
$$
$$ 
 \triangledown ^2 F=
\left(
\begin{array}{cc}
	6x^2+2y^2-2 & 4xy \\
	4xy& 2x^2+6y^2+2
\end{array}
\right) \approx
\left(
\begin{array}{cc}
	-2 & 4xy \\
	4xy& 2
\end{array}
\right).
$$

The eigenvalues of $ \triangledown ^2 F$ are
$$\lambda_1=2(2x^2-\sqrt{x^4+2x^2y^2-2x^2+y^4+2y^2+1}+2y^2) \approx -2$$ 
and 
$$\lambda_2=2(2x^2+\sqrt{x^4+2x^2y^2-2x^2+y^4+2y^2+1}+2y^2) \approx 2.$$ 
The eigenvectors of $ \triangledown ^2 F$ are
$$
u_1=
\left(
\begin{array}{cc}
	-1-\sqrt{1+4x^2y^2} \\
	2xy
\end{array}
\right) \approx
\left(
\begin{array}{cc}
	-1 \\
	0
\end{array}
\right) 
$$ 
and 
$$
u_2=
\left(
\begin{array}{cc}
	2xy \\
	1+\sqrt{1+4x^2y^2}
\end{array}
\right) \approx
\left(
\begin{array}{cc}
	0 \\
	1
\end{array}
\right).
$$
Note that $\|\nabla F(z_0)\|\approx 0$ as long as $x, y$ are sufficiently small.
Then, \begin{equation*}
\begin{aligned}
        A_0:&=\nabla^2 F(z_0)+\delta_j \|\nabla F(z_0)\|^{\tau}Id \approx \nabla^2 F(z_0),
\end{aligned}
\end{equation*}
and 
\begin{equation*}
\begin{aligned}
v_0:&= A_0^{-1}\triangledown F(z_0)\approx
\left(
\begin{array}{cc}
	-\frac{1}{2} & 0\\
  0 & \frac{1}{2}
\end{array}
\right)
 \cdot
\left(
\begin{array}{cc}
	-2x\\
	2y
\end{array}
\right)\\
&=
\left(
\begin{array}{cc}
	x \\
	y
\end{array}
\right)=-x \cdot u_1+y \cdot u_2.
\end{aligned}
\end{equation*}
Since the eigenvalue with respect (approximately) to $u_1$ is negative, and the eigenvalue with respect (approximately) to $u_2$ is positive, we have 
$$w_0=x \cdot u_1+y \cdot u_2=(-x,y). 
$$

By definition of BNQN New Variant, $z_1=z_0-\gamma w_0$ for some $\gamma >0$. Therefore, $x_1>x_0$ as wanted. 

In conclusion, if we start any point $z_0\in H_1=\{(x,y):~x>0\}$, the sequence $\{z_n\}$ constructed by BNQN New Variant will converge to the root $z_1^*=1$. 

Since the function $F(x,y)$ is invariant under the rotation $R(x,y)=(-x,y)$, which maps $H_1$ to $H_2=\{(x,y):~x<0\}$ and maps $z_1^*$ to $z_2^*=-1$, we conclude that if $z_0\in H_2$ then the sequence $\{z_n\}$ constructed by BNQN New Variant will converge to the root $z_2^*=-1$. This completes the proof of part 2a of Theorem \ref{TheoremMain1}. 

b) Now, we show that if $z_0\in L=\{(x,y):~x=0\}$, then the sequence $\{z_n\}$ constructed by BNQN New Variant will converge to the point $0$. 

Indeed, as in the proof of Case 2 in part 2a of Theorem \ref{TheoremMain1}, we can show that in this case $z_n\in L$ for all $n$. We also know that $\{z_n\}$ must converge to either $-1$, $0$ or $1$. Since, only the point $0$ belongs to $L$, it follows that $\{z_n\}$ must converge to $0$. 
\end{proof}

{\bf Remark:} Inspired by the proofs of Subcase 2.1 and Claim B, one may wonder if it is true that for any initial point $z_0=(x,y)$ with $0<x<<1$ and $y>0$, and for the point $z_1=(x_1,y_1)$ constructed by BNQN New Variant from the point $z_0$, then we should always have $x_1\geq x_0$. However, this speculation is not true, as one can easily check with computer experiments.

\end{document}